\documentclass[]{article}
\usepackage{amsmath}
\usepackage{amssymb}

\newcommand{\be}{\begin{equation}}
\newcommand{\ee}{\end{equation}}

\newtheorem{definition and theorem}[theorem]{Definition and Theorem}

\newtheorem{*remark}[theorem]{$^* $Remark}

\newtheorem{*exercise}[theorem]{$^* $Exercise}
\newtheorem{**exercise}[theorem]{$^{** } $Exercise}

\begin{document}
\title{Approximation of Multiple Integrals over Hyperboloids with
Application to a Quadratic Portfolio with Options. }

\author{Jules Sadefo Kamdem
\thanks{
Paper based on the M. Jules SADEFO KAMDEM University of Reims Phd Thesis. }\\
Universite De Reims  \\
Laboratoire de Math\'ematique\\ UMR 6056-CNRS\\
UFR Sciences Reims  \\
BP 1039 Moulin de la Housse\\ 51687 Reims Cedex FRANCE\\
sadefo@univ-reims.fr\\
\and
Alan Genz\\
Department of Mathematics\\
Washington State University \\
Pullman, WA 99164-3113 USA\\
alangenz@wsu.edu}

\maketitle
\begin{abstract}
We consider an application involving a financial quadratic
portfolio of options, when the joint underlying log-returns
changes with multivariate elliptic distribution. This motivates
the need for methods for the approximation of multiple integrals
over hyperboloids. A transformation is used to reduce the
hyperboloid integrals to a product of two radial integrals and two
spherical surface integrals. Numerical approximation methods for
the transformed integrals are constructed. The application of
these methods is demonstrated using some financial applications
examples.
\end{abstract}

\section{Introduction}

Value-at-Risk (VaR) is considered to be one of the standard measures of market
risk. VaR measures the maximum loss that a portfolio can experience
with a certain probability over a certain horizon, for example, one day.
Mathematically, if the profit or loss is given by $\Pi(t)-\Pi(0)$,
for which $\Pi(t)$ is the price of the portfolio at $t$, VaR for a
confidence level 1-$\alpha $ is determined by the following equation:
$$
  \mathbb{P}rob\{\Pi(0,S(0))-\Pi(t,S(t))>VaR_{\alpha} \}=\alpha
  \label{prob},
$$
where  $S(t)=( S_{1}(t), \ldots ,S_{n}(t))$ is a vector of asset prices that
govern risk factors.

In this paper, we consider numerical methods for the estimation of
integrals over hyperboloids and their application to
VaR computations for a complex portfolio that contains options depending on
the market fluctuations that create risk. We reduce the problem of VaR
computations to multiple integrals over hyperboloids, and show how these
integrals can be approximated using techniques described by
Genz and Monahan \cite{GM} and Sheil and O'Muircheartaigh\cite{SO}.

One of the most important analytic methods for VaR computation, which is
called  $\Delta$-normal VaR, was introduced in the
RiskMetrics Technical Document (1996). The method is based on the assumptions
that the distribution is Normal and the portfolio is linear.
Sadefo-Kamdem \cite{S2}, generalized the $\Delta$-Normal
VaR by introducing the $\Delta$-elliptic VaR for a linear portfolio,
with the $\Delta$-Student VaR given as an example.
An advantage of the $\Delta$-elliptic VaR (for example, $\Delta$-Normal
or $\Delta$-Student) is that the formula is still fairly simple to calculate.
But in practice, if we deal with a $\Delta$-hedged portfolio,
the $\Delta$-elliptic VaR does not provide a realistic model, and that is why
alternatives were proposed in the paper of Brummelhuis, Cordoba,
Seco and Quintanilla \cite{BCQS}. In that paper,
mathematical stationary phase techniques were used
to estimate VaR for a quadratic portfolio when the risk factor changes
with a Normal distribution. In a sequel paper, Brummelhuis and
Sadefo-Kamdem \cite{BS} provided an estimation method with more
precision for a VaR with a quadratic portfolio and
generalized Laplace distributions.

We assume, as in \cite{BS}, that the approximation for the price of the
portfolio is given by
$$
\Pi(t)-\Pi(0)\approx \Theta t + \Delta  \mathbb{X}^t
      + \frac{1}{2}\mathbb{X}^t \Gamma \mathbb{X},
$$
and we also assume that the joint log-returns $\mathbb{X}$ is
elliptically distributed.
For further details about elliptic distributions, see
Embrechts, McNeil and Straussman \cite{EMS}.
$\Gamma$, $\Theta$ and $\Delta$ are functions of some sensitivities
of the portfolio (see Taleb \cite{TA}, 1997, for a discussion concerning
sensitivities).
We also suppose that $t=1$, because the time horizon for VaR is
generally taken to be one day
If the log-returns of $\mathbb{X}$ are Normally distributed,
Albanese and Seco \cite{AS} have shown how to reduce
the analysis of a quadratic VaR to the computation of the integral of
a Gaussian over a quadric in a space of possibly very high dimension.
We will use the more general assumption of an elliptic
distribution for the risk factors.

This paper  proposes numerical methods for the approximation of integrals over
hyperboloids, with application to estimate the VaR. We combine some techniques
described in \cite{G} ,\cite{GM} and \cite{SO} to approximate the integrals
over hyperboloids for VaR, with the generalized assumption that the underlying
joint log-returns changes with an elliptic distribution.
To illustrate our method, we will take the familiar case of Normal
distributions, and we consider test examples for two
$\Delta$-hedged portfolios from the French CAC 40 market.
Brummelhuis, Cordoba, Quintanilla and Seco \cite{BCQS}
have considered a quadratic portfolio with an analytic
approximation for the Gaussian integrals over quadrics.
Sadefo-Kamdem and Brummelhuis \cite{BS} have provided a similar analysis
with a Generalized Laplace distribution. Albanese and Seco \cite{AS}
investigated the approximation of integral of a Gaussian
over a hyperboloid region with Fourier transform methods.

One of the most common methods for quadratic
perturbations of the linear VaR uses the
Cornish-Fisher expansion for the quantile function of
non-Gaussian variables. There are also some quadratic
approximations in Hull \cite{HU} and Dowd \cite{D}.
Many papers in literature have proposed numerical
methods for the quadratic approximation (see, for example, \cite{S1},
where Sadefo-Kamdem proposed the use of some numerical methods of
Genz, \cite{G}, and the use of hypergeometric functions for
a portfolios of equities VaR with multivariate {t-Student distribution}).

The rest of our paper is organized as follows. In Section
2 and 3, following Albanese and Seco \cite{AS}, we
show how portfolio volatility  can be used to reduce the
calculation of VaR to the approximation of integral over
hyperboloid, assuming elliptic distributions that admit a density function.
In section 4, we propose a numerical method for the  approximation of
integrals over hyperboloids  using some methods of Genz \cite{G},
Genz and Monahan \cite{GM} and \cite{SO}.
To illustrate our method we use examples where the density function is
Normal. In section 5, we consider two examples of financial
portfolios, and we have showed that our method is applicable to
estimate the VaR for the portfolio. In
Section 6 , we provide some conclusions.

\section{Quadratic Portfolio of Options Application}

In this section, we will define a quadratic portfolio of options as
Quintanilla did in \cite{QU}. We first define
$\mathbb{X}=(\mathbb{X}_{1},\ldots,\mathbb{X}_{n})$, with
$$ \mathbb{X}_{i}=log(S_{i}(t)/S_{i}(0)), $$
and we define
$\Delta_{1}= (\Delta_{1}^1 , \ldots, \Delta_{1}^n )$, with
$\Delta_{1}^i = S_{i}(0).\Delta^i $, and
$\Delta =(\Delta^1 , \ldots , \Delta^n ) $, the gradient vector of the
portfolio at time $t=0$.
We also define
$\Gamma_1=\Big{(} \Gamma_1^{i,j}\Big{)}_{i,j=1,...,n}$ by
$$
  \Gamma_1^{i,j} =\left \{
   \begin{array}{ll}
  S_{i}^2 (0) \Gamma^{i,i} + \Delta_1^i  & \textrm{ if $i=j$ }\\
  S_{i} (0) S_{j} (0) \cdot \Gamma^{i,j} & \textrm{ if $i \neq j$  }
  \end{array} \right. ,
$$
with
$\Gamma=\Big{(} \Gamma^{i,j}\Big{)}_{i,j=1,...,n}=
\Big{(} \frac{\partial^2 \Pi }{\partial S_i \partial S_j }(0)
\Big{)}_{i,j=1,...,n}$, the  Hessian of the portfolio at time $t=0$.
If we use a $2^{nd}$ order Taylor series approximation for $\Pi$, then
$$
\Pi(t,S(t))-\Pi(0,S(0))\approx t\Theta + \Delta_{1} \mathbb{X}^t
+ \mathbb{X}\Gamma_{1} \mathbb{X}^{t}/2 \label{eq0},
$$
where $\Theta=\frac{\partial \Pi}{\partial t}(0) \label{b1}$ .

If we consider a $\Delta$-hedged Portfolio, we have
$\Delta=0$, and therefore $\Delta_{1}=0$.
Our goal is to determinate the Value-at-Risk quantity $V$ with confidence
level $1-\alpha$ when $t=1$, as a solution to the
equation
$$
G_{\Gamma_1}(-V)= \mathbb{P}(\Theta + \mathbb{X} \Gamma_{1}
\mathbb{X}^{t}/2 \leq -V)=\alpha \label{eq1}.
$$

If we assume that the joint underlying log-returns $\mathbb{X}$ have a
multivariate elliptic distribution with zero mean,
then $G_{\Gamma_{1}}(-V)$ is given by
$$
G_{\Gamma_1}(-V)= \int_{\{\Theta + x\Gamma_{1}
x^t /2 \leq -V \}} g(x \Sigma^{-1}x^t) \frac{dx}{\sqrt{det(\Sigma)}}=\alpha .
$$

\section{Transformation to a Hyperboloid Integration Region}

We first decompose $\Sigma$ as $\Sigma=\mathbb{C}\mathbb{C}^t$, where
$\mathbb{C}$ is the (lower-triangular) Cholesky decomposition
factor of $\Sigma$, and then we use the transformation $x =
y\mathbb{C}^t$ to give
$$
G_{\Gamma}(-V)=\int_{\{\Theta+y\mathbb{C}^t\Gamma\mathbb{C}y^{t}/2\leq
-V \}} g(yy^t) dy .
$$
We next assume that the sensitivity-adjusted variance-covariance
matrix, $\mathbb{C}^t\Gamma\mathbb{C}$, has a diagonalization in
the form $\mathbb{C}^t\Gamma\mathbb{C}=
\mathbb{O}\mathbb{D}\mathbb{O}^t$, with $\mathbb{O}$ orthogonal,
and $\mathbb{D}$ diagonal. $\mathbb{C}\Gamma \mathbb{C}^{t}$ is
not necessarily definite, but we can assume that a diagonalization
has been constructed with
$$
\mathbf{\mathbb{D}}=
\left(\begin{array}{ccc}D_+ & 0 \\ 0  &  -D_-\end{array} \right),
$$
where
$$
D_{\epsilon}= \left(\begin{array}{ccc}
d_{1}^{\epsilon} & \ldots & 0\\
 0    &\ddots& 0 \\
 0 & \ldots& d_{n_{\epsilon}}^{\epsilon}
\end{array} \right)
$$
for $\epsilon=\pm{1}$, and where all $d_{j}^{+}$,$d_{j}^{-} \geq
0$, and $-d_{1}^{-}\leq -d_{2}^{-}\leq \ldots \leq
-d_{n_{-}}^{-}\leq d_{1}^{+}\leq \ldots \leq  d_{n_{+}}^{+}$ .
Then, we can use the transformation $z =y\mathbb{O}$ to give
$$
G_{\Gamma}(-V)=\int_{\{\Theta+z\mathbb{D}z^{t}/2\leq -V \}}g(zz^t)
dz ,
$$
and finally, we can use the transformation $w =|\mathbb{D}|^{1/2}z$ to give
\begin{equation}
G_{\Gamma}(-V)= \int_{\{|w_{+}|^2 -|w_{-}|^2 \leq -2(V + \Theta )
\}} g( w|\mathbb{\mathbb{D}}|^{-1}w^{t})
\frac{dw}{\sqrt{det(|\mathbb{D}|)}}, \label{eq3}
\end{equation}
where $w=(w_{+},w_{-})$ is the decomposition of $\mathbb{R}^n $
into the respective positive and negative subspaces of the
eigenbasis for $\mathbb{C}^t\Gamma\mathbb{C}$. After changing the
direction of the inequality, we obtain the following expression for
$G(R)$, which will be the starting point for the discussion of our
computational methods.
\begin{equation}
G(R)=   \int_{\{|w_{-}|^{2} - |w_{+}|^{2}  \geq  R^{2}  \} } g(w
|\mathbb{D}|^{-1}w^{t})
  \frac{dw}{\sqrt{det(|\mathbb{D}|)}},    \label{eq4}
\end{equation}
where $R^2 = 2(V+\Theta)$. The integration region is the {\it hyperboloid}
defined by $|w_{-}|^2 - |w_{+}|^2  \geq  R^{2}$.
Our goal is to determine $R$ a solution to $G(R)=\alpha$. Once
we find $R$, we will have the approximate quadratic Value-at-Risk given by
$V = R^2/2 - \Theta$.

\section{Integration over Hyperboloids }
\subsection{The Normal Case}
For many applications, the distribution $g$ is a Normal
distribution. In these cases,
$$
G(R)=\int_{\{|w_{-}|^{2} - |w_{+}|^{2}  \geq  R^{2}  \} } e^{-w
|\mathbb{D}|^{-1}w^{t}/2} \frac{dw}{\sqrt{(2\pi)^n
det(|\mathbb{D}|)}} .
$$
Separating the $w$ variables, we find
$$
G(R)= \int_{\{|w_+|^2  \geq  0  \} } e^{-w_+ D_+^{-1}w_+^{t}/2}
\int_{\{|w_-|^2   \geq  R^{2} + |w_{+}|^{2} \} } e^{-w_-
D_-^{-1}w_-^{t}/2} \frac{dw_-}{\sqrt{(2\pi)^{n_{-}} det(D_-)}}
\frac{dw_+}{\sqrt{(2\pi)^{n_+} det(D_+)}} .
$$
The inner integral for $w_-$ can be efficiently computed using the
algorithm described by Sheil and O'Muircheartaigh \cite{SO}, so we
define $H(R,r)$ by
$$
H(R,r) = \int_{\{|w_-|^2   \geq  R^{2} + r^2 \} } e^{-w_-
D_-^{-1}w_-^{t}/2} \frac{dw_-}{\sqrt{(2\pi)^n_{-} det(D_-)}},
$$
and let $w_+ = D_+^{\frac{1}{2}}z$
Then $G(R)$ can rewritten as
$$
G(R) = \int_{-\infty}^\infty \int_{-\infty}^\infty \ldots\int_{-\infty}^\infty
e^{-xx^t/2} H(R,|zD_+z^t|)\frac{dz} {\sqrt{(2\pi)^{n_+}}} .
$$
Integrals in this form can be approximated using methods described by
Genz and Monahan \cite{GM}.

\subsection{The General Case}
We need to determine approximations to integrals in the form
$$
G(R) = \int_{\{|x|^{2} - |y|^{2} \geq R^2 \}}\varphi(x,y) dy dz,
$$
where $x \in \mathbb{R}^{n_1}$ and $y \in \mathbb{R}^{n_2}$. If we
use the changes of variable: $y = r_2\xi_2$, $x = r_1\xi_1$, with
$r_2=|y|$ and $r_1=|x|$, then $G(R)$ becomes
$$
G(R) = \int_{0}^{\infty}r_2^{n_2-1} \int_{|\xi_2|=1}
\int_{\sqrt{R^2+r_2^2}}^\infty r_1^{n_1-1}\int_{|\xi_1|=1}
\varphi(r_1\xi_1,r_2\xi_2) d\sigma(\xi_1) dr_1 d\sigma(\xi_2) dr_2
$$
We now have $G(R)$ defined in terms of a product of two integrals
over hyper-spherical surfaces, defined by $|\xi_2|=1$ and
$|\xi_1|=1$, and two radial integrals. The hyper-sphere surface
integrals can be approximated using methods described in the paper
by Genz \cite{G}. If the surface and radial integrals are combined,
then generalizations of the methods described by Genz and Monahan \cite{GM}
can be used. Efficient approximation of the radial integrals
will depend on information about the rate of decrease of the
integrand $\varphi$ for large values of $r_1$ and $r_2$.

\section{Application Examples}
We will distinguish 3 case in our analysis :
\begin{itemize}
\item $n_{-}$=0; if $g$ is Normal,
$$
G(R)= \int_{\{|w_+|^2   \leq  R^2\} } e^{-w_+D_+^{-1}w_+^{t}/2}
\frac{dw_+}{\sqrt{(2\pi)^{n_{+}} det(D_+)}}.
$$
$G(R)$ can be efficiently computed using the Sheil and O'Muircheartaigh
\cite{SO} algorithm, when $R^2 = -2(V+\Theta) \geq 0$.
\item $n_{+}$=0; if $g$ is Normal,
$$
G(R)= \int_{\{|w_-|^2   \geq  R^{2}  \} } e^{-w_-
D_-^{-1}w_-^{t}/2} \frac{dw_-}{\sqrt{(2\pi)^{n_{-}} det(D_-)}}.
$$
$G(R)$ can be efficiently computed using the Sheil and O'Muircheartaigh
\cite{SO} algorithm.
\item $n_{-}$ and $n_{+}$ are both nonzero; $g$ is Normal,
$$
G(R) = \int_{-\infty}^\infty \int_{-\infty}^\infty \ldots\int_{-\infty}^\infty
e^{-xx^t/2} H(R,|zD_+z^t|)\frac{dz} {\sqrt{(2\pi)^{n_+}}} .
$$
Integrals in this form can be approximated using methods that are a
combination of the  Sheil and O'Muircheartaigh \cite{SO} and
Genz and Monahan \cite{GM} algorithms.
\end{itemize}
\subsection{An Example when $n_{+}$=0 }
We construct a $\Delta$-hedged portfolio that contains $n$ equities and
$n$ European call options on these equities from the
{\it French CAC 40 Market}. The Price of the Portfolio is given  by
 $$\label{b0}
\Pi(t,S(t))=\sum_{i=1}^n [ -C_{i}(t,S_{i}(t))+ \Delta^i \cdot S_{i}] ,
$$
where $S_{i}$ is an equity price $i$, with $S(t)=(S_{1},\ldots,S_{n})$,
and $C_{i}(t,S_{i}(t))$ is the price of European call option $i$ on
equity $i$. $\Delta$ is known in the literature as a
gradient portfolio sensitivity vector. Our portfolio has been chosen so
that $\Delta=0$,
with $\Delta^i=\frac{\partial C_{i}}{\partial S_{i}}(S_{i}(0))$, and
$\Delta=(\Delta^1,\ldots,\Delta^{n})$.
The exercise price of each European call option is given in the following
table for an example where $n = 9$.:
\bigskip
\begin{center}
Table 1: Data for Nine CAC 40 European Call Options
    \begin{tabular}{|c|c|c|c|c|}
      \hline
       & Exercise Price & Interest Rate & Maturity & Underlying Price\\
      \hline
      BNPPARIBAS  & 44.26 &0.1 & 3 month& 39.75\\
      \hline
      BOUYGUES & 23.49&0.1 & 3 month&27.30 \\
      \hline
      CAP GEMINI & 34.71 &0.1 & 3 month&24.00 \\

      \hline
      CREDIT AGRICOLE & 17.36&0.1 & 3 month&14.80 \\
      \hline
      DEXIA & 11.5 &0.1 & 3 month&9.38\\
      \hline
      LOREALL & 61.85 &0.1 & 3 month&62.90 \\
      \hline
      TF1 & 26.38 &0.1 & 3 month&22.02\\
      \hline
      THOMSON & 15.22 & 0.1 & 3 month&17.13\\
      \hline
      VIVENDI & 16.19 &0.1 & 3 month&17.00 \\
      \hline
    \end{tabular}
\end{center}
\bigskip
Using the above data with the exponential moving
weighted average (EMWA), we obtained the following $\Sigma$:
$$
\Sigma= \left( \begin{array}{rrrrrrrrr}
 0.0017 &-0.0001 &0.0012 & 0.0005 &0.0008 &0.0008 &0.0008 &0.0002 & 0.0002 \\
-0.0001 & 0.0009 &0.0005 &-0.0001 &0.0000 &0.0000 &0.0000 &0.0006 & 0.0005 \\
 0.0012 & 0.0005 &0.0038 & 0.0006 &0.0011 &0.0008 &0.0014 &0.0006 & 0.0007 \\
 0.0005 &-0.0001 &0.0006 & 0.0006 &0.0002 &0.0004 &0.0004 &0.0001 & 0.0001 \\
 0.0008 & 0.0000 &0.0011 & 0.0002 &0.0015 &0.0007 &0.0008 &0.0000 & 0.0002 \\
 0.0008 & 0.0000 &0.0008 & 0.0004 &0.0007 &0.0011 &0.0006 &0.0004 & 0.0005 \\
 0.0008 & 0.0000 &0.0014 & 0.0004 &0.0008 &0.0006 &0.0013 &0.0000 & 0.0001 \\
 0.0002 & 0.0006 &0.0006 & 0.0001 &0.0000 &0.0004 &0.0000 &0.0019 & 0.0007 \\
 0.0002 & 0.0005 &0.0007 & 0.0001 &0.0002 &0.0005 &0.0001 &0.0007 & 0.0029 \\
\end{array} \right)
$$
The matrix $\Gamma$ is a diagonal matrix with diagonal entries given by
$$
d=(-116.4889, -33.4063, -11.8389, -21.1723, -11.9582, -161.2178,
-34.7884, -19.7664, -27.2993)
$$
and $ \Theta=-31.2689. $
The eigenvalues of the $\mathbb{D}$ matrix are given by the vector
$$e=(-.05025,-.1456,-.0605,-.0424,-.0231,-.0.0053,-.0101,-.0154,-.0131).$$
The following Table provides some $R$ and $V$ values
that were found as numerical solutions to the equation $G(R)=\alpha$,
for selected $\alpha$'s.

\begin{center}
\begin{tabular}{|c|c|c|c|}\hline
 $\alpha$&   0.05& 0.025 &  0.01 \\  \hline
 $R$     & 0.9160&1.0038&1.1128\\ \hline
 $V$     & 31.6883& 31.7727&31.8881\\  \hline
\end{tabular}
\end{center}

\subsection{Example with $n_{-}>0$ and $ n_{+}>0$  }
We consider a portfolio that contains call options and
   put options on equities, so that price of the portfolio at time t,
   is given by:
   $$ \Pi(t)=\sum_{i=1}^5 [C_{i}(t,S_{i}(t)) -\delta_{i} S_{i}(t)] -
    \sum_{j=6}^{10} P_{j}(t,S_{j}(t)) - (\delta_{j}-1) S_{j}(t)]  $$
The prices of each of the options are taken from data from the
{\it French CAC 40 Market}, and are given in the following table:
\bigskip
\begin{center}
Table 2: Data for Ten CAC 40 European Call and Put Options
    \begin{tabular}{|c|c|c|c|c|}
      \hline
       & Exercize Price & Interest Rate & Maturity & Underlying Price\\
      \hline
     Call-BNPPARIBAS  & 30.00 &0.05 & 3 months & 39.75\\      \hline
     Call-BOUYGUES & 19.00 &0.05 & 3 months & 27.30\\      \hline
     Call-CAP GEMINI & 20.00 & 0.05 & 3 months & 24.00 \\      \hline
     Call-CREDIT AGRICOLE & 10.50 &0.05 & 3 months & 14.80\\      \hline
      Call-DEXIA & 9.00 &0.05 & 3 months & 9.38\\      \hline
      Call-LOREALL & 40.00 &0.05 & 3 months & 62.90 \\      \hline
      Put-SOCIETEGENERALE & 50& 0.05 &3 months & 64.00\\      \hline
      Put-TF1 & 18.00 &0.05 & 3 months & 22.02 \\      \hline
      Put-THOMSON & 9.00 & 0.05 & 3 months  & 17.13\\      \hline
      Put-VIVENDI &9.00 &0.05 & 3 months & 17.00\\      \hline
    \end{tabular}
\end{center}
\bigskip
In this example, using the three month historical
data for the ten CAC 40 equities, with the exponential moving weighted average
(EMWA) and $\lambda=0.94$, we obtained the following $\Sigma$.
$$
\Sigma= \left( \begin{array}{rrrrrrrrrr}
 0.0016&-0.0001&0.0012& 0.0005& 0.0008& 0.0007&-0.0000& 0.0008&0.0002& 0.0002\\
-0.0001& 0.0008&0.0004&-0.0001& 0.0000& 0.0000& 0.0004&-0.0000&0.0005& 0.0005\\
 0.0011& 0.0004&0.0035& 0.0006& 0.0010& 0.0007& 0.0004& 0.0013&0.0006& 0.0006\\
 0.0005&-0.0001&0.0006& 0.0005& 0.0002& 0.0004&-0.0001& 0.0003&0.0001& 0.0001\\
 0.0008& 0.0000&0.0010& 0.0002& 0.0015& 0.0007&-0.0000& 0.0008&0.0000& 0.0002\\
 0.0007& 0.0000&0.0007& 0.0004& 0.0007& 0.0010&-0.0002& 0.0005&0.0003& 0.0004\\
-0.0000& 0.0004&0.0004&-0.0001&-0.0000&-0.0002& 0.0015&-0.0002&0.0008&-0.0001\\
 0.0008&-0.0000&0.0013& 0.0003& 0.0008& 0.0005&-0.0002& 0.0012&0.0000& 0.0001\\
 0.0002& 0.0005&0.0006& 0.0001& 0.0000& 0.0003& 0.0008& 0.0000&0.0018& 0.0007\\
 0.0002& 0.0005&0.0006& 0.0001& 0.0002& 0.0004&-0.0001& 0.0001&0.0007& 0.0026\\
\end{array} \right).
$$
The matrix $\Gamma$ is a diagonal matrix with diagonal
$$
d=(24.186, 8.6269, 21.7320, 4.3111, 15.4949,
-4.5815,-82.2915, -22.2079, -1.2957, -1.2822),
$$
and $\Theta=-3.8596.$
The eigenvalues of the $\mathbb{D}$ matrix are given by the vector
$$
e=(-0.1251,-0.0115,-0.0030,-0.0014,-0.0006,0.0014,0.0092,0.0124,0.0290,0.1271).
$$
The following Table provides some $R$ and $V$ values
that were found as numerical solutions to the equation $G(R)=\alpha$,
for selected $\alpha$'s.

\begin{center}
\begin{tabular}{|c|c|c|c|}\hline
 $\alpha$&   0.05& 0.025 &  0.01 \\  \hline
 $R$     & 0.6069&0.7176&0.8455\\ \hline
 $V$     & 4.0438&4.1171&4.2166\\  \hline
\end{tabular}
\end{center}

\section{Conclusions}
We have considered an application in the domain of Risk Management
to estimate the Value-at-Risk of the portfolio of an option. We have
reduced the problem to a multiple integral over a hyperboloid. This
type of integral can be approximated  using techniques described by
Genz and Monahan \cite{GM} and Sheil and O'Muircheartaigh \cite{SO}.

\end{document}